\documentclass[12pt]{amsart}
\usepackage[utf8]{inputenc}
\usepackage{amsmath,amssymb,amsthm}

\newtheorem{theorem}{Theorem}
\theoremstyle{definition}
\newtheorem{example}{Example}
\newtheorem{question}{Question}

\title{A note on a sumset in $\mathbb{Z}_{2k}$}

\author{
Octavio A. Agustín-Aquino
}
\address{Universidad Tecnológica de la Mixteca\\
Instituto de Física y Matemáticas\\
Carretera a Acatlima Km. 2.5\\
Huajuapan de León, Oaxaca, México, C.P. 69000}
\email{octavioalberto@mixteco.utm.mx}

\date{15 December 2017}

\begin{document}

\begin{abstract}
Let $A$ and $B$ be additive sets of $\mathbb{Z}_{2k}$, where $A$ has cardinality $k$ and 
$B=v.\complement A$ with $v\in\mathbb{Z}_{2k}^{\times}$. In this note some bounds for the
cardinality of $A+B$ are obtained, using four different approaches. We also prove that in a special case
the bound is not sharp and we can recover the whole group as a sumset.
\end{abstract}
\maketitle

\section{Introduction}

In music, a \emph{canon} is typically understood as a musical composition where a melody is
imitated by various voices, with a duration offset between them (well known examples are ``Row, Row, Row Your Boat'' or ``Frère Jacques'').
A canon like those can be more aptly described as a ``pitch canon'', in contraposition to the ``rhythmic canons'' introduced by Oliver Messiaen,
where the rhythm is imitated instead of the melody. In this direction, a remarkable and pioneering use of sumsets in music was done by Dan Tudor Vuza, introducing what he
called ``Regular Complementary Canons of Maximal Category'', which are aperiodic sumsets $S,T\subseteq \mathbb{Z}_{n}$
such that $\mathbb{Z}_{n} = S+T$, where $S$ (or $T$) represents the set of duration offsets between rhythmic imitations.
For a nice introduction to the fascinating interplay of music and mathematics
in this regard, see \cite{mA11}.

In this note, we make further connections between sumsets and the musical realm of \emph{counterpoint}, where canon is
but one of its techniques. Thus let $U$ and $V$ be additive subsets of $\mathbb{Z}_{2k}$ with cardinality $k$, and
\begin{align*}
 U+V&=\{u+v:u\in U,v\in V\},\\
 x.U&=\{xa:a\in U\}.
\end{align*}

I stumbled upon the
problem of proving that, if $k$ is large enough and under certain hypothesis regarding the structure of $U$, we have
\[
U+V = \mathbb{Z}_{2k}
\]
where $U$ is a set closely related to $V$. Hence $U$ and $V$ are akin to Vuza canons, except for aperiodicity is not
required beforehand but some other conditions are to be fulfilled. To be more specific, a very interesting case from the
mathematical counterpoint theory perspective is when
\[
V = v.\complement U, \quad v\in\mathbb{Z}_{2k}^{\times}\setminus\{-1\}
\]
and, additionally, $\complement U = U+k$ (here $\complement$ stands for the set complement with
respect to $\mathbb{Z}_{2k}$). In order
to explain why, let $\overrightarrow{GL}(\mathbb{Z}_{2k})$ be the set of bijective functions
\begin{align*}
 e^{u}.v:\mathbb{Z}_{2k}&\to\mathbb{Z}_{2k},\\
 x&\mapsto vx+u,
\end{align*}
where $v\in\mathbb{Z}_{2k}^{\times}$ and $u\in\mathbb{Z}_{2k}$. If $A\subseteq \mathbb{Z}_{2k}$ is such
that $g(A)\neq A$ for every $g\in \overrightarrow{GL}(\mathbb{Z}_{2k})$ except the identity, and $A\cup p(A)=\mathbb{Z}_{2k}$ for a unique $p\in \overrightarrow{GL}(\mathbb{Z}_{2k})$, then it
is called a \emph{counterpoint dichotomy} and $p$ is its \emph{polarity}.

\begin{example} One interesting example is $A=\{0,2,3\}\subseteq \mathbb{Z}_{6}$, whose
polarity is $e^{1}.-1$, since it is essentially the only one available. Another important specimen is
\[
 K =\{0,3,4,7,8,9\}\subseteq \mathbb{Z}_{12},
\]
with polarity $e^{2}.5$, for $K$ is the set of consonances in Renaissance counterpoint modulo
octave, when the intervals in $12$-tone equally tempered scale are interpreted as $\mathbb{Z}_{12}$.
The interested reader may consult \cite[Part VII]{mazzola2002topos} or \cite{AJM15} and references therein for further details.
\end{example}

Throughout this paper, we will attack (with varying degrees of generality) the following question.

\begin{question}\label{Q:Importante}
Given a subset $A\subseteq \mathbb{Z}_{2k}$ of cardinality $k$, is it true that
\begin{equation}\label{E:CondF}
 A+v.\complement A =\begin{cases}
 \mathbb{Z}_{2k}, &v\in\mathbb{Z}_{2k}^{\times}\setminus \{-1\},\\
 \mathbb{Z}_{2k}\setminus\{0\}, &v=-1?
 \end{cases}
\end{equation}
\end{question}

When this question can be answered in the affirmative then, for any $e^{u}.(-v)$ except the identity, 
there exists $x\in A$ and $y\in \complement A$ such that
\[
x+(-v)y = u\quad\text{or}\quad vy+u = x\quad\text{or}\quad e^{u}.(-v)(y)\in  A
\]
which means that no element of $\overrightarrow{GL}(\mathbb{Z}_{2k})$ but the identity leaves the set $A$ invariant.
If there exists also a $p\in\overrightarrow{GL}(\mathbb{Z}_{2k})$ such that $p(A)=\complement A$,
then $A$ is a counterpoint dichotomy.

A set that I have been trying to prove is a counterpoint dichotomy for a long time (some
reasons for this are stated in \cite{oA17})
via answering Question \ref{Q:Importante} is
\begin{equation}\label{E:IntSet}
A = \{0,1\}\cup \{3,4,\ldots,k-1\}\cup \{k+2\}.
\end{equation}

It is not difficult to verify that $e^{k}.1(A)=\complement A$ and to see that
\[
A+A = \mathbb{Z}_{2k}\quad\text{and}\quad A-A \supseteq \mathbb{Z}_{2k}\setminus\{k\},
\]
since $2=1+1$, $2k-1 = (k+2)+(k-3)$ and $k+1=(k-2)+3$ for the first equality. The other
one is consequence of $3-1=2$ and $1-3 = -2$.

Although the following three sections do not prove $A$ satisfies the rest of \eqref{E:CondF},
they provide some evidence and results that may be interesting on their own. Moreover,
an elementary proof of this fact found by Merlijn Staps is presented in the fifth section. In
the last section, some final remarks are made.

\section{Using the Ruzsa distance}
Let $U$ and $V$ be subsets of an additive group $G$. A couple of weak bounds for $|U+V|$ can be obtained using Ruzsa's useful notion of 
``distance'' in additive combinatorics
\[
d(U,V) = \log\frac{|U-V|}{\sqrt{|U||V|}},
\]
which is a seminorm. In particular, it satisfies a triangle inequality
\[
d(U,V) \leq d(U,W)+d(W,V).
\]

Note now that, regarding the set \eqref{E:IntSet}, we have
\[
d(A,-A)= \log\frac{|A+A|}{|A|} = \log\frac{2k}{k} = \log 2;
\]
the number $\delta(U)= \exp(d(U,-U))$ is the \emph{doubling constant} of the set $U$, and thus $\delta(A)=2$.

From the Ruzsa triangle inequality we can deduce \cite[p. 61]{TV06}
\[
|U||V-V|\leq |U+V|^{2}
\]
which, for the case of $V=A$ and $U=B$, specializes to
\[
|A+B| \geq \sqrt{|B||A-A|}\geq \sqrt{k(2k-1)} = \sqrt{2-\frac{1}{k}}k.
\]

On the other hand, again by the triangle inequality
\[
\log 2 = d(A,-A) \leq d(A,B)+d(B,-A)
\]
and a pigeon-hole argument, either
\[
d(A,B)\geq \frac{1}{2}\log 2
\]
or
\[
d(-A,B)=d(A,-B)\geq \frac{1}{2}\log 2.
\]

Equivalently, either
\[
|A-B|\geq \sqrt{2}k
\]
or
\[
|A+B|\geq \sqrt{2}k.
\]

We conclude that, for any subsets $A$ and $B$ of the cardinality $k$ such that $\delta(A)=2$, we have
\[
\max\{|A+B|,|A-B|\}\geq \sqrt{2}k.
\]

We do not know if there exist pairs of subsets of $\mathbb{Z}_{2k}$ such that $A$ has doubling
constant $2$ and $|A+B|$ or $|A-B|$ get
arbitrarily close to this bound.

\section{Using additive energy and a theorem by Olson}

Let
\[
[P]=\begin{cases}
1, & P\text{ is true},\\
0, & \text{otherwise},
\end{cases}
\]
be the Iverson bracket \cite[p. 24]{GKP94}, and define
the \emph{additive energy} of the subsets $U$ and $V$ of the additive group $G$
by
\[
E(U,V) = \sum_{u_{1},u_{2}\in U,v_{3},v_{4}\in V}[u_{1}+u_{2}=v_{3}+v_{4}].
\]

Another well-known inequality \cite[p. 63]{TV06} for the cardinality of $U+V$ is
\[
|U\pm V| \geq \frac{(|U||V|)^{2}}{E(U,V)}.
\]

From this we infer another strategy to improve the previous estimates for $|A+B|$, namely finding upper
bounds for $E(A,B)$. A good start might be the Cauchy-Schwarz inequality
\[
E(A,B)\leq \sqrt{E(A,A)E(B,B)}.
\]

This seems promising when $B=v.\complement A$ and $\complement A = A+\{k\}$, since the invertibility of $v$ implies
\begin{align*}
E(v.\complement A,v.\complement A)&=\sum_{a_{1},a_{2},a_{3},a_{4}\in A+\{k\}}[va_{1}+va_{2}=va_{3}+va_{4}]\\
&=\sum_{a_{1},a_{2},a_{3},a_{4}\in A}[v(a_{1}+a_{2})=v(a_{3}+a_{4})]\\
&=\sum_{a_{1},a_{2},a_{3},a_{4}\in A}[a_{1}+a_{2}=v^{-1}v(a_{3}+a_{4})]\\
&=\sum_{a_{1},a_{2},a_{3},a_{4}\in A}[a_{1}+a_{2}=a_{3}+a_{4}]=E(A,A).
\end{align*}

Thus $E(A,v.\complement A)\leq E(A,A)$. Nevertheless, this straightforward approach loses some of its charm as soon as we calculate a few values
of the energy and the corresponding bounds.

\begin{table}
\begin{center}
\begin{tabular}{|c|c|c|c|}
\hline
$k$ & $E(A,A)$ & $\frac{(|A||v.\complement A|)^{2}}{E(A,A)}=\frac{k^4}{E(A,A)}$ & $\frac{k^{3}}{2E(A,A)}$\\
\hline
$8$ & $296$ & $13.84$ & $0.86$\\
\hline
$9$ & $425$ & $15.44$ & $0.86$\\
\hline
$10$ & $590$ & $16.95$& $0.85$\\
\hline
$11$ & $795$ & $18.42$& $0.84$\\
\hline
$12$ & $1044$ & $19.86$& $0.83$\\
\hline
$100$ & $665180$ & $150.34$& $0.751$\\
\hline
$1000$ & $666651080$ & $1500.04$& $0.750$\\
\hline
\end{tabular}
\caption{Additive energy $E(A,A)$ for small $k=|A|$, where $A$ is defined by \eqref{E:IntSet}, the
corresponding bound for $|A+B|$ and the fraction of $\mathbb{Z}_{2k}$ that is thus guaranteed to be covered by $A+B$.}
\label{T:Calc}
\end{center}
\end{table}

As it is readily seen in Table \ref{T:Calc}, the quality of the bound is expected to decrease as $k$ increases, although it would remain as a mild improvement
with respect the one obtained in the previous section. In fact, assuming $E(A,A)$ is a polynomial in $k$, from a simple interpolation from the
data in Table \ref{T:Calc} we find that
\[
 E(A,A) = \tfrac{2}{3}k^{3}-\tfrac{47}{3}k+80.
\]

This means that, for $k\geq 6$, we have $E(A,A)\leq \frac{2}{3}k^{3}$, and then
\[
 |A\pm v.\complement A| \geq \tfrac{3}{2}k.
\]

This bound can be obtained from a theorem due to Olson, and actually it holds for any set $B$ of cardinality $k$, not only
those of the form $v.\complement A$. Before stating 
Olson's theorem, observe that an additive subset $U$ of $G$
is contained in a coset of a unique smallest subgroup $H$ of $G$. Denote with $[U]$ such a coset.
\begin{theorem}[Olson, 1984, \cite{olson1984sum}, \cite{vsevolod2006critical},\cite{boothby2013new}]
Let $U$ and $V$ be additive subsets of $G$. If $U+V\neq G$ and $[U]=G$, then $|U+V|\geq \tfrac{1}{2}|U|+|V|$.
\end{theorem}

Suppose $G=\mathbb{Z}_{2k}$ and $U=A$. Any coset containing $A$ has cardinality at least $k$. But it cannot
have exactly $k$ elements, for the cosets would be forced to be either the set of even elements of $\mathbb{Z}_{2k}$
or its complement, but clearly $A$ is contained in neither. Thus $[A]=\mathbb{Z}_{2k}$, so if $A+B$
is not the whole group, it must consist of at least $\tfrac{1}{2}k+k=\tfrac{3}{2}k$ elements.

\section{Using trigonometric sums}
Let $r_{U+V}(t)$ the number of representations of $t$ as a sum $t=u+v$ for $u\in U$ and
$v\in V$, where $U$ and $V$ are additive subsets of a group $G$. The following is a standard
technique using the so-called trigonometric sums in number theory (a readable and short introduction can be found in \cite{MG09}). Note first that
\[
\frac{1}{m}\sum_{\xi=0}^{m-1}e^{2\pi i \xi x /m} = [x\equiv 0 \pmod{m}],
\]
so we can write
\[
\frac{1}{2k}\sum_{\xi=0}^{2k-1}e^{2\pi i \xi(u+v-\lambda)/(2k)}=[u+v\equiv \lambda\pmod{2k}].
\]

If we sum over $U$ and $V$ and exchange the order of summation,
\begin{align*}
r_{U+V}(\lambda) &=\sum_{u\in U}\sum_{v\in V}[u+v\equiv \lambda\pmod{2k}].\\
&= \frac{1}{2k}\sum_{u\in U}\sum_{v\in V}\sum_{\xi=0}^{2k-1}e^{2\pi i \xi(u+v-\lambda)/(2k)}\\
&= \frac{1}{2k}\sum_{\xi=0}^{2k-1}\left(\sum_{u\in U}e^{2\pi i \xi u/(2k)}\sum_{v\in V}e^{2\pi i \xi v/(2k)}\right) e^{-2\pi i \xi \lambda/(2k)},
\end{align*}
and then we extract the $\xi=0$ term, we conclude
\[
r_{U+V}(\lambda) =  \frac{k}{2}+E
\]
where, by the triangle inequality,
\begin{equation}\label{E:Error}
|E|\leq \frac{1}{2k}\sum_{\xi=1}^{2k-1}\left|\sum_{u\in U}e^{\pi i \xi u/k}\right|\left|\sum_{v\in V}e^{\pi i \xi v/k}\right|=2k\sum_{\xi=1}^{2k-1}|\widehat{1}_{U}(\xi)||\widehat{1}_{V}(\xi)|
\end{equation}
and
\[
 \widehat{f}(\xi):=\frac{1}{|G|}\sum_{x\in G}f(x)\overline{e^{2\pi i \xi x/|G|}}
\]
is the Fourier transform. Observe now that $|\widehat{1}_{v.\complement A}(\xi)|=|\widehat{1}_{\complement A}(\xi)|\leq |\widehat{1}_{A}(\xi)|$, so
for $U=A$ and $V=v.\complement A$, we have
\[
|E|\leq 2k\sum_{\xi=1}^{2k-1}|\widehat{1}_{A}(\xi)|^2 \leq k,
\]
which is not useful. On the other hand, since (see \cite[Lemma 6, Chapter 1]{vinogradov1954method})
\[
|\widehat{1_{A}}(\xi)| \leq \frac{1}{2k\sin(\pi\xi/(2k))}+\frac{1}{k}
\]
then
\begin{align*}
\sum_{\xi=1}^{2k-1}|\widehat{1}_{A}(\xi)|^{2}&\leq\sum_{\xi=1}^{2k-1}\left(\frac{1}{2k\sin(\pi\xi/(2k))}+\frac{1}{k}\right)^{2}\\
&= 2\sum_{\xi=1}^{k}\left(\frac{1}{2k\sin(\pi\xi/(2k))}+\frac{1}{k}\right)^{2}-\frac{9}{4k^{2}}.
\end{align*}

Now the sequence
\[
a_{k,\xi} = \begin{cases}
\left(\frac{1}{2k\sin(\pi\xi/(2k))}+\frac{1}{k}\right)^{2}, &1\leq\xi\leq k,\\
0,&\text{otherwise},
\end{cases}
\]
is such that $a_{k,\xi}\geq a_{k+1,\xi}$ and $\sum_{\xi=1}^{\infty}a_{1,\xi}=\frac{9}{4}$. By the monotone convergence theorem,
we obtain
\[
\lim_{k\to\infty}\sum_{\xi=1}^{2k-1}|\widehat{1}_{A}(\xi)|^2 = 2\lim_{k\to \infty}\sum_{\xi=1}^{k-1}\frac{1}{\pi^{2}\xi^{2}} = \frac{1}{3},
\]
which amounts to estimate $|E|\leq \tfrac{2}{3}k$ for large $k$, but that is not enough to ensure that $r_{A+v.\complement A}(\lambda)\geq 0$ for any $\lambda$ and $v\neq -1$.
Furthermore, it suggests that the most we can get this way is $|A+v.\complement A|\geq \frac{5}{6}k$ (see \cite[p. 210]{TV06}).

\section{Using a result by Mann}

For a penultimate attempt we use the following generalization of the celebrated Cauchy-Davenport theorem.

\begin{theorem}[Mann, 1965, see \cite{SZ00}]
Let $S$ be a subset of an arbitrary abelian group $G$. Then one of the following holds:
\begin{enumerate}
\item For every subset $T$ such that $S+T \neq G$, we have $|S+T|\geq |S|+|T|-1$.
\item There exists a proper subgroup $H$ of $G$ such that $|S+H|< |S|+|H|-1$.
\end{enumerate}
\end{theorem}

Thus one of these two alternatives holds:

\begin{enumerate}
\item It is true that $|A+v.\complement A|\geq |A|+|v.\complement A|-1 = 2k-1$.
\item There is proper subgroup $H$, such that
\[
 |A+H|< k+|H|-1.
\]
\end{enumerate}

We claim that, for the set $A$, we have
\[
|A+v.\complement A|\geq 2k-1
\]
by discarding the second alternative. In order to do so, suppose $H=\langle d\rangle$ where $0\leq d\leq k$ and
\[
|H+A|< k+|H|-1.
\]

Being that $H$ is proper, we have $|H|\leq k$. Let us suppose that $d\geq 1$ (since
the trivial case is evidently false), which implies that $|H|=\frac{2k}{d}$. Thus $A+H$ is
the placement of copies of $A$ with spaces of $d$ elements, so it covers all the elements of $\mathbb{Z}_{2k}$
with at most $\frac{2k}{d}$ exceptions, thus
\[
 k+\frac{2k}{d}-1> |A+H|\geq 2k-\frac{2k}{d}.
\]

This is possible if, and only if,
\[
 \frac{2k}{d}+k-1 > 2k-\frac{2k}{d}
\]
or, equivalently,
\[
 4 > \frac{4k}{k+1} > d,
\]
thereof $d=2$ or $d=3$. If $d=2$, we are done, for $A$ has $\{0,1\}$ as a subset, thus $A+H=\mathbb{Z}_{2k}$,
a contradiction.

In the later case (which arises only when $3$ divides $k$), it would be possible that each ``slot'' of $3$ elements
$\{3j,3j+1,3j+2\}$ determined by $H$ and to be covered by $A$ to have
$3j+2$ uncovered. Nevertheless, the ``antipodal'' slot $\{3j+k,3j+k+1,3j+k+2\}$ would not allow this to happen, since the
potentially uncovered element must be covered with the translate
$(3j+k+2)+k$ of $3j+k+2\in A+3j$. Moreover: we
are certain that a copy of $A$ is placed in $k$ because $3$ is one of its factors. So, $A+H$ would leave no
element uncovered, for there are an even number of slots, each one paired with its antipode.
Hence $H=\langle 3\rangle$ is also an impossibility.

From the above proof we also obtain that $A$ is aperiodic, i. e., $A+H\neq A$ except for $H=\{0\}$.
Invoking Kemperman structure theorem (as stated, for example, in \cite[p. 71-72]{grynkiewicz2009step})\footnote{More
specifically, the pair $(A,-v.\complement A)$ is of type IV in the classification stated in \cite[p. 71]{grynkiewicz2009step}.}, we conclude that
\[
 A-\complement A = \mathbb{Z}_{2k}\setminus\{0\}
\]
and, furthermore, if $A+v.\complement A\neq \mathbb{Z}_{2k}$, then there exists $u$ such that
\[
 v.\complement A = u-\complement A.
\]

This equivalent to the following: except for for $v=-1$,
and $u=0$ it is true that
\[
 -v.A+u \neq A,
\]
which means exactly that $A$ is a counterpoint dichotomy. Thus, Kemperman's theorem cannot lead us further in relation to
the cardinality of $A+v.\complement A$.

\section{A proof for a special case}

Question \ref{Q:Importante} can be answered affirmatively for the set \eqref{E:IntSet}, when
$k\geq 10$, as we now show. Let us first identify $v$ with an element in
\[
 \{-k+1,-k+2,\ldots,-1,0,1,\ldots,k-1\}.
\]

Observe that for $v=-1$ we have $A+(-1).B=\mathbb{Z}_{2k}\setminus \{k\}$, and for $v=1$
we have $A+1.B=A+B=A+k+A=\mathbb{Z}_{2k}+k =\mathbb{Z}_{2k}$.
Now suppose $k-3\geq |v|>1$, that is, $k-3\geq |v|\geq 3$. Choose
\[
 X=\{3,4,\ldots,k-1\}
\]
and $Y=X-3$. We claim that $Y+v.Y=\mathbb{Z}_{2k}$, for this would imply that
\begin{align*}
 A+v.A&\supset X+v.X\\
&= Y+3+v.(Y+3)\\
&=Y+v.Y+3+3v\\
&= \mathbb{Z}_{2k}+3+3v =\mathbb{Z}_{2k},
\end{align*}
and hence $A+v.B = A+v.(A+k) = A+v.A+vk = \mathbb{Z}_{2k}$, as we want.

To prove the claim, we note that the set $Y+v.Y$ contains the multiples $0,v,2v,\ldots,(k-4)v$. We have
\[
 |(k-4)v|>|3(k-4)|\geq 2k
\]
and, since $|v|\leq k-3$, the elements between multiples of $v$ are also in $Y+v.Y$, as $Y$ contains $\{0,1,\ldots,v-1\}$.
This means that $Y+v.Y=\mathbb{Z}_{2k}$. The remaining cases we need to deal with are
\[
 v\in \{\pm(k-1),\pm(k-2)\};
\]
note that $\pm(k-2)$ only occurs when $k$ is odd.

For $v=k-1$, we note that $A+(k-1).A$ contains $A$ and $A+(k-1)=A+k-1=B-1$, so $A+(k-1).A$ contains all the elements
of $\mathbb{Z}_{2k}$ with the possible exceptions of those in $B\setminus(B-1)$. However
\begin{align*}
-1 &= (k+2)+(k-1)3,\\
 k+1 &= 4+(k-1)3,\\
 2 &= 6+(k-1)4,
\end{align*}
proving that all of them belong to $A+(k-1).A$. We must have $4,6\in A$ since $k\geq 10$.

For $v=-(k-1) = k+1$ we have $B\setminus (B+1) = \{2,k,k+2\}$, and the analogous certificates are
\begin{align*}
2 &= (k-1)+(k+1)3,\\
k &= (k-4)+(k+1)4,\\
k+2 &= (k-3)+(k+1)4;
\end{align*}
for $v=k-2$ we have $B\setminus (B-2) = \{-1,-2,k,2\}$ and
\begin{align*}
-1 &= 1+(k-2)1,\\
-2 &= 0+(k-2)1,\\
k &= 6+(k-2)3,\\
2 &= 10+(k-2)4;
\end{align*}
finally, for $v=k+2$ we have $B\setminus (B+2) = \{2,k,k+1,k+4\}$ and
\begin{align*}
2 &= k+(k+2)1,\\
k &= k+(k+2)0,\\
k+1 &= (k-7)+(k+2)4,\\
k+4 &= (k-4)+(k+2)4.
\end{align*}

\section{Some final remarks}

The results distilled from Mann's and Kemperman's theorems take us rather close to the goal of
proving that \eqref{E:CondF} holds for the set $A$ defined by \eqref{E:IntSet}, but ultimately fail.
We can manage to provide an elementary proof of the fact, but we do not know how much this approach
can be generalized, or what this means for the classificatory nature of Kemperman's theorem.

Nevertheless, these facts make evident that there is a significant gap between $E(A,v.\complement A)$
and $E(A,A)$. They also point out that, in order to succeed with the use of exponential sums, a very sharp
estimate of \eqref{E:Error} is required.
 
\section{Acknowledgments}

I am indebted to Merlijn Staps at Universiteit Utrecht for his kind help reading an earlier version
of this paper and finding a proof of the sumset property of \eqref{E:IntSet}.

I sincerely thank my friends José Hernández Santiago at Centro de Ciencias Matem\'aticas (UNAM) and 
Marcelino Ramírez Ibáñez at Universidad del Papaloapan for their invaluable feedback regarding early
drafts of this paper. I also thank an anonymous reviewer who brought an embarrassing blunder to my attention,
as well as another who helped me to polish my prose and to clarify several obscure points.

\bibliographystyle{amsplain}
\bibliography{simple}

\end{document}